\input amstex
\magnification=\magstep1 
\baselineskip=13pt
\documentstyle{amsppt}
\vsize=8.7truein \CenteredTagsOnSplits \NoRunningHeads

\def\PP{\Cal P}
\def\QQ{\Cal Q}
\def\UU{\frak U}
 \topmatter
 
\title Computing the partition function for graph homomorphisms with multiplicities \endtitle 
\author Alexander Barvinok and Pablo Sober\'on \endauthor
\address Department of Mathematics, University of Michigan, Ann Arbor,
MI 48109-1043, USA \endaddress
\email barvinok$\@$umich.edu, psoberon$\@$umich.edu  \endemail
\date July 2015 \enddate
\thanks  The research of the first author was partially supported by NSF Grants DMS 0856640 and DMS 1361541.
\endthanks 
\abstract We consider a refinement of the partition function of graph homomorphisms and present a quasi-polynomial algorithm to compute it in a certain domain. As a corollary, we obtain quasi-polynomial algorithms for computing partition functions for independent sets, perfect matchings, Hamiltonian cycles and dense subgraphs in graphs as well as for graph colorings. This allows us to tell apart in quasi-polynomial time graphs that are sufficiently far from having a structure of a given type (i.e., independent set of a given size, Hamiltonian cycle, etc.) from graphs that have sufficiently many structures of that type, even when the probability to hit such a structure at random is exponentially small.

\endabstract

\endtopmatter

\document

\head 1. Introduction and main results \endhead

\subhead (1.1) Partition function of graph homomorphisms with prescribed multiplicities \endsubhead Let $G=(V, E)$ be an undirected graph with set $V$ of vertices, set $E$ of edges, without loops or multiple edges. We denote by $\Delta(G)$ the largest degree of a vertex in $G$. In what follows, we assume that $\Delta(G) \geq 1$, so that the graph contains at least one edge.

Let $m=\left(\mu_1, \ldots, \mu_k \right)$ be a vector of positive integers such that 
$$\mu_1 + \ldots + \mu_k = |V|$$ and let 
$A=\left(a_{ij}\right)$ be a $k \times k$ symmetric real or complex matrix. We define the 
{\it partition function of graph homomorphisms of $G$ with multiplicities $m$} by 
$$\PP_{G,m}(A)=\sum \Sb \phi:\ V \rightarrow \{1, \ldots, k\} \\ \left|\phi^{-1}(i)\right|=\mu_i \text{\ for \ } i=1, \ldots, k \endSb
\prod_{\{u, v\} \in E} a_{\phi(u) \phi(v)}. \tag1.1.1$$
Here the sum is taken over all maps $\phi: V \longrightarrow \{1, \ldots, k\}$ that map precisely $\mu_i$ vertices of $G$ into $i$ for 
all $i=1, \ldots, k$ and the product is taken over all edges of $G$. If we take the sum over all $k^{|V|}$ maps $\phi: V \longrightarrow \{1, \ldots, k\}$, without the multiplicity restrictions, we obtain what is known as the {\it partition function of graph homomorphisms}
$$\PP_G(A) =\sum_{\phi:\ V \rightarrow \{1, \ldots, k\}} \prod_{\{u, v\} \in E} a_{\phi(u) \phi(v)}, \tag1.1.2$$
see \cite{BG05} and \cite{C+13}.
\bigskip
$\bullet$ One of the main results of our paper is a deterministic algorithm, which, given a graph $G=(V, E)$, a vector 
$m=\left(\mu_1, \ldots, \mu_k\right)$ of multiplicities, a $k \times k$ symmetric matrix $A=\left(a_{ij}\right)$ such that 
$$\left| a_{ij} - 1 \right| \ \leq \ {\gamma \over \Delta(G)} \quad \text{for all} \quad i, j \tag1.1.3$$
and a real $\epsilon >0$, computes $\PP_{G, m}(A)$ within relative error $\epsilon$ in $\left(|E|k\right)^{O(\ln |E| -\ln \epsilon)}$ time.
Here $\gamma >0$ is an absolute constant (we can choose $\gamma=0.1$). More precisely, we prove that for {\it complex} matrices $A$ satisfying (1.1.3) we necessarily have $\PP_{G, m}(A) \ne 0$ and, given a real $\epsilon >0$, we compute in $\left(|E|k\right)^{O(\ln |E| -\ln \epsilon)}$ time a (possibly complex) number $\tilde{\PP}_{G,m}(A) \ne 0$ such that 
$$\left| \ln \PP_{G, m}(A) - \ln \tilde{\PP}_{G,m}(A)\right| \ < \ \epsilon$$
(we choose a branch of the logarithm which is real when $a_{ij}=1$ for all $i,j$).
\bigskip
In \cite{BS15}, we construct a deterministic algorithm, which, given a graph $G=(V, E)$ and a $k \times k$ symmetric matrix 
$A=\left(a_{ij}\right)$ satisfying (1.1.3) with a better constant $\gamma=0.34$ (we can choose $\gamma=0.45$ if $\Delta(G) \geq 3$ and $\gamma=0.54$ for all sufficiently large $\Delta(G)$) computes the graph homomorphism partition 
function (1.1.2) within relative error $\epsilon$ in $\left(|E|k\right)^{O(\ln |E| -\ln \epsilon)}$ time. Although the methods of \cite{BS15} and this paper are similar, it appears that neither result follows from the other. 

Specializing (1.1.1), we obtain various quantities of combinatorial interest. 

\subhead (1.2) Independent sets in graphs \endsubhead Let $G=(V, E)$ be a graph as above. A set $S \subset V$ is called {\it independent} if $\{u, v\} \notin E$ for every two vertices $u, v \in S$. Finding the maximum size of an independent set in a given graph within a factor of $|V|^{\epsilon}$ is an NP-hard problem for any $0 \leq \epsilon < 1$, fixed in advance \cite{H\aa 99}, \cite{Zu07}. 

Let us choose $k=2$ and define the matrix $A=\left(a_{ij}\right)$ by 
$$a_{11}=0, \quad a_{12}=a_{21}=a_{22}=1.$$
Let $m=\left(\mu_1, \mu_2\right)$ be an integer vector such that $\mu_1 + \mu_2 =|V|$. One can see that a map 
$\phi: V \longrightarrow \{1, 2\}$ contributes $1$ to (1.1.1) if $\phi^{-1}(1)$ is an independent set in $G$ and contributes 0 otherwise. Hence $\PP_{G,m}(A)$ is the number of independent sets in $G$ of size $\mu_1$. 

Let us modify $A \mapsto \tilde{A}$ by choosing 
$$\tilde{a}_{11} = 1-{\gamma \over \Delta(G)} \quad \text{and} \quad \tilde{a}_{12}=\tilde{a}_{21}=\tilde{a}_{22} =1 +{\gamma \over \Delta(G)},$$
where $\gamma$ is the constant in (1.1.3). Then 
$$\left(1+{\gamma \over \Delta(G)}\right)^{-|E|} \PP_{G,m}(\tilde{A}) = \sum \Sb S: \ S \subset V \\ |S| =\mu_1 \endSb w(S), \tag1.2.1$$
where 
$$w(S) =\left(1 +{\gamma \over \Delta(G)}\right)^{-e(S)} \left(1 -{\gamma \over \Delta(G)}\right)^{e(S)}$$
and $e(S)$ is the number of pairs of vertices of $S$ that span an edge of $G$. 
Thus 
$$\exp\left\{ -2\gamma {e(S) \over \Delta(G)} - \gamma^3 {e(S) \over \Delta^3(G)} \right\} \ \leq \ w(S) \ \leq \ \exp\left\{-2 \gamma {e(S) \over \Delta(G)}\right\},$$
so roughly
$$w(S) \approx \exp\left\{ -2 \gamma {e(S) \over \Delta(G)}\right\}.$$
Hence (1.2.1) computes a weighted sum over all subsets $S$ of vertices of cardinality $\mu_1$, where $w(S)=1$ if $S$ is an independent set in $G$, while all other subsets are weighted down exponentially in the number of edges that the vertices of the subsets span. 

Computing (1.2.1) allows us to distinguish graphs with sufficiently many independent sets of a given size from graphs that are sufficiently far from having an independent set of a given size. Indeed, if every subset $S \subset V$ of $\mu_1$ vertices of $G$ spans at least $x$ edges of $G$, we get 
$$\left(1+{\gamma \over \Delta(G)}\right)^{-|E|} \PP_{G,m}(\tilde{A}) \ \leq \ {|V| \choose \mu_1}
 \exp\left\{-2\gamma {x \over \Delta(G)}\right\}. $$
If, on the other hand, a random subset of $\mu_1$ vertices of $G$ is an independent set with probability at least 
$\displaystyle 2\exp\left\{-2\gamma x \over \Delta(G) \right\}$ then 
$$\left(1+{\gamma \over \Delta(G)}\right)^{-|E|} \PP_{G,m}(\tilde{A}) \ \geq \ 2{|V| \choose \mu_1} \exp\left\{-2\gamma {x \over \Delta(G)}\right\}.$$
Computing $\PP_{G,m}(\tilde{A})$ within a relative error of $0.1$, say, we can distinguish these two cases in $|E|^{O(\ln |E|)}$ time.

Let us fix some $\delta > 0$ and $\epsilon_0 >0$ and let us consider a class of graphs $G$ that satisfy $|E| \geq \delta |V| \Delta(G)$ (in other words, we consider graphs that are not very far from regular). 
Let us choose $x = \epsilon |E|$ for some $\epsilon \geq \epsilon_0 >0$. 
In this case, $G$ is sufficiently far from having an independent set of size $\mu_1$ if every subset of vertices of size $\mu_1$ spans at least some constant proportion $\epsilon \geq \epsilon_0$ (arbitrarily small, but fixed in advance) of the total number $|E|$ of edges whereas $G$ has sufficiently many independent sets if the probability that a randomly selected set of $\mu_1$ vertices is independent is at least 
$$ 2 \exp\left\{ -2 \gamma \epsilon {|E| \over \Delta(G)}\right\}  \ \leq \ 2 \exp\left\{ -2 \gamma \epsilon_0 \delta |V| \right\}.$$
We can distinguish these two cases in quasi-polynomial time, although in the latter case the probability to hit an independent set at random is exponentially small in the number of vertices of the graph.

More generally, we can separate graphs where each subset $S \subset V$ of cardinality $\mu_1$ spans at least $x$ edges of $G$ from graphs having sufficiently many subsets $S \subset V$ of cardinality $\mu_1$ spanning at most $y < x$ edges of $G$.

\subhead (1.3) Hafnians, Hamiltonian permanents, and subgraph densities \endsubhead

Suppose that $G$ is a union of $n$ pairwise vertex-disjoint edges, so $|V|=2n$ and $\Delta(G)=1$. Suppose further that $k=|V|$ and let us choose 
$m=(1, \ldots, 1)$. Then, up to a normalizing constant, $\PP_{G,m}(A)$ is the {\it hafnian} of the matrix $A$, see, for example, 
Section 8.2 of \cite{Mi78}. In particular, if $A$ is the adjacency matrix of a simple undirected graph $H$ with set $\{1, \ldots, k\}$ of vertices, then the value of 
$(2^n n!)^{-1} \PP_{G,m}(A)$ is the number of {\it perfect matchings} in $H$, that is, the number of collections of edges of $H$ covering every vertex of $H$ exactly once. Hence we obtain a deterministic algorithm approximating the hafnian of a $2n\times 2n$ symmetric matrix $A=\left(a_{ij}\right)$ that satisfies 
$$\left| 1- a_{ij} \right| \ \leq \ \gamma \quad \text{for all} \quad i, j \tag1.3.1$$
within a relative error $\epsilon >0$ in $n^{O(\ln n - \ln \epsilon)}$ time. A similar algorithm, though with a better constant 
$\gamma=0.19$ (which can be improved to $0.27$, see Lemma 4.1 below and the remark thereafter), was earlier constructed in \cite{B15a}.

Suppose that $G$ is a cycle with $n>2$ vertices, so that $\Delta(G)=2$, and that $k=n$. Let us choose $m=(1, \ldots, 1)$. If $A$ is the adjacency matrix of a simple undirected graph $H$ then the value of $(2n)^{-1}\PP_{G,m}(A)$ is the number of {\it Hamiltonian cycles} in $H$, that is, the number of walks that visit every vertex in $H$ exactly once before returning to the starting point. For a general $A$, the value of $n^{-1} \PP_{G,m}(A)$ is a Hamiltonian version of the permanent of $A$, 
$${1 \over n} \PP_{G,m}(A)=\sum_{\sigma} \prod_{i=1}^n a_{i \sigma(i)},$$
where the sum is taken over all $(n-1)!$ permutations $\sigma$ of the set $\{1, \ldots, n\}$ that consist of a single cycle,
cf. \cite{Ba15}.
 Hence we obtain a deterministic algorithm approximating the Hamiltonian permanent of an $n \times n$ symmetric matrix $A=\left(a_{ij}\right)$ that satisfies 
$$\left| 1- a_{ij} \right| \ \leq \ {\gamma \over 2} \quad \text{for all} \quad i, j $$
within a relative error $\epsilon >0$ in $n^{O(\ln n - \ln \epsilon)}$ time. This result is new, although there is a polynomial time algorithm, which, for, given a real $A$ satisfying (1.3.1) with any $\gamma <1$, fixed in advance, approximates the Hamiltonian permanent of (not necessarily symmetric) $A$ within a factor of $n^{O(\ln n)}$ 
(the implicit constant in the ``$O$" notation depends on $\gamma$), see \cite{Ba15}. As in Section 1.2, we obtain a polynomial time algorithm that distinguishes graphs that have sufficiently many (for example, at least $\epsilon^n n!$ for some $0 < \epsilon <1$, fixed in advance) Hamiltonian cycles from graphs that are sufficiently far from Hamiltonian (need at least $\epsilon n$ new edges added to become Hamiltonian for some fixed $0 < \epsilon <1$), see \cite{Ba15} for details.

For an arbitrary $G$ with $n$ vertices, let $k \geq n$ and let $A$ be the adjacency matrix of a simple graph $H$ with set $\{1, \ldots, k\}$ of vertices. Let us consider a graph $\widehat{G}$ consisting of $G$ and $k-n$ isolated vertices. Let us choose a $k$-dimensional vector $m=\left(1, \ldots, 1\right)$. 
Then $\bigl( (k-n)!\bigr)^{-1} \PP_{\widehat{G}, m}(A)$ is the number of embeddings $\psi: G \longrightarrow H$ that map distinct vertices of $G$ into distinct vertices of $H$ and the edges of $G$ into edges of $H$. The case of a complete graph $G$ is of a particular interest. Assuming that $G$ is a complete graph with $n$ vertices and ${n \choose 2}$ edges, we conclude that 
 $\bigl( n! (k-n)!\bigr)^{-1} \PP_{\widehat{G}, m}(A)$ is the number of {\it cliques} of size $n$ in $H$. In this case we have $\Delta(\widehat{G})=n-1$.  Given $H$ with vertex set $\{1, \ldots, k\}$, let us modify $A \longmapsto \tilde{A}$ by 
$$\tilde{a}_{ij}=\cases 1+ {\gamma \over n-1} &\text{if $\{i, j\}$ is an edge of $H$} \\ 1-{\gamma \over n-1} &\text{if $\{i, j\}$ is not an edge of $H$.} \endcases$$
Given a set $S$ of vertices of $H$, let $t(S)$ be the number of pairs of distinct vertices of $S$ that do not span an edge of $H$. Then
$$(n!(k-n)!)^{-1} \left(1+{\gamma \over n-1}\right)^{-{n \choose 2}} \PP_{\widehat{G}, m}(\tilde{A}) =
\sum \Sb S: \ S \subset \{1, \ldots, k\} \\ |S|=n \endSb w(S), \tag1.3.2$$
where 
$$w(S)=\left(1+{\gamma \over n-1}\right)^{-t(S)} \left(1-{\gamma \over n-1}\right)^{t(S)}=
\exp\left\{- 2\gamma {t(S)  \over n-1} + O\left({t(S) \over (n-1)^3}\right)\right\}.$$
Hence we obtain a deterministic algorithm of $(kn)^{O(\ln n)}$ complexity, which, given a graph $H$ with $k$ vertices, computes (within a relative error of $0.1$, say) the sum (1.3.2) of the weights of $n$-subsets $S$ of the set of vertices, where each weight $w(S)$ is exponentially small in the number of edges of $H$ that the vertices of $S$ fail to span. A similar algorithm was constructed earlier in \cite{B15b} (in \cite{B15b} a worse constant $\gamma=0.07$ is achieved in the general case; however,  in the most interesting case of $n=o(k)$ and $n \geq 10$,  \cite{B15b} achieves a better constant of $\gamma=0.27$). Again, it allows us to distinguish graphs with no dense induced subgraphs of size $n$ from graphs having sufficiently many dense induced subgraphs of size $n$, even when ``many" still allows for the probability to hit a dense induced subgraph at random to be exponentially small in $n$, see \cite{B15b}.
Finding the densest induced subgraph of a given size $n$ in a given graph with $k$ vertices is a notoriously hard problem. The best known algorithm of $k^{O(\ln k)}$ complexity approximates the highest density of an $n$-subgraph up to within a multiplicative factor of $k^{1/4}$ \cite{B+10}, \cite{Bh12}.

Although an independent set in a graph corresponds to a clique in the complementary graph on the same set of vertices, there appears to be no direct relation between the partition functions (1.2.1) and (1.3.2). In (1.2.1), the weight of a subset depends on the maximum degree $\Delta(G)$ while in (1.3.2) it depends on the size of the subset.

\subhead (1.4) Graph colorings \endsubhead Let $G=(V, E)$ be a graph as in Section 1.1 and let 
$$\phi: V \longrightarrow \{1, \ldots, k\}$$
be a map.
The map $\phi$ is called a {\it proper $k$-coloring} if $\phi(u) \ne \phi(v)$ whenever $\{u, v\} \in E$. The smallest $k$ for which a proper $k$-coloring of $G$ exists is called the {\it chromatic number} of $G$. Approximating the chromatic number of a given graph within a factor of $|V|^{1-\epsilon}$ for any $0 < \epsilon < 1$, fixed in advance, is NP-hard \cite{FK98}, \cite{Zu07}.
Let us define a matrix $A=\left(a_{ij}\right)$ by 
$$a_{ij}=\cases 1 &\text{if $i\ne j$} \\ 0 &\text{if\ } i=j.\endcases$$ Given an integer vector $m=\left(\mu_1, \ldots, \mu_k \right)$, we observe that $\PP_{G, m}(A)$ is the number of proper colorings of $G$, where the $i$-th color is used exactly $\mu_i$ times. Colorings with prescribed number of vertices of a given color were studied in the case of {\it equitable colorings}, where 
any two multiplicities $\mu_i$ and $\mu_j$ differ by at most 1, see \cite{K+10} and references wherein.

As above, let us modify $A\longmapsto \tilde{A}$ by letting 
$$\tilde{a}_{ij}=\cases 1+{\gamma \over \Delta(G)} &\text{if\ } i \ne j \\ 1-{\gamma \over \Delta(G)} &\text{if\ } i=j. \endcases$$
Then we obtain
$$\left(1+{\gamma \over \Delta(G)}\right)^{-|E|} \PP_{G,m}(\tilde{A})= \sum \Sb \phi: \ V \rightarrow \{1, \ldots, k\} \\ \left| \phi^{-1}(i)\right|=\mu_i 
\ \text{for} \ i=1, \ldots, k \endSb w(\phi),$$
where 
$$w(\phi)=\left(1+{\gamma \over \Delta(G)}\right)^{-e(\phi)} \left(1-{\gamma \over \Delta(G)} \right)^{e(\phi)} \tag1.4.1$$
and $e(\phi)$ is the number of miscolored edges of $G$ (that is, edges whose endpoints are colored with the same color under the coloring $\phi$).
Thus 
$$\exp\left\{ -2\gamma {e(\phi) \over \Delta(G)} - \gamma^3 {e(\phi) \over \Delta^3(G)} \right\} \ \leq \ w(e) \ \leq \ \exp\left\{-2 \gamma {e(\phi) \over \Delta(G)}\right\}$$
and hence (1.4.1) represents a weighted sum over all colorings into $k$ colors with prescribed multiplicity of each color and the weight of each coloring being exponentially small in the number of miscolored edges. As before, we can compute (1.4.1) within a relative error of $0.1$ in $(k|E|)^{O(\ln |E|)}$ time.

\subhead (1.5) Partition function of  edge-colored graph homomorphisms with multiplicities \endsubhead
It turns out that instead of the partition function $\PP_{G,m}(A)$ defined by (1.1.1), it is more convenient to consider a more general expression. Let $G=(V, E)$ be a graph as in Section 1.1 and let $B=\left(b^{uv}_{ij}\right)$ be $|E| \times {k(k+1) \over 2}$ real or complex matrix with entries indexed by edges $\{u, v\} \in E$ and unordered pairs $1 \leq i, j \leq k$. Technically, we should have written $b^{\{u, v\}}_{\{i, j\}}$, but we write just $b^{uv}_{ij}$, assuming that 
$$b^{uv}_{ij}=b^{vu}_{ij}=b^{vu}_{ji}=b^{uv}_{ji}.$$
We define
$$\QQ_{G,m}(B)=\sum \Sb \phi:\ V \rightarrow \{1, \ldots, k\} \\ \left|\phi^{-1}(i)\right|=\mu_i \text{\ for \ } i=1, \ldots, k \endSb
\prod_{\{u, v\} \in E} b_{\phi(u) \phi(v)}^{uv}. \tag1.5.1$$
Let $H$ be a simple undirected graph with set $\{1, \ldots, k\}$ of vertices and suppose that the edges of $G$ and $H$ are colored. Let us define 
$$b^{uv}_{ij}=\cases 1 &\text{if\ } \{u, v\} \ \text{and} \ \{i, j\} \ \text{are edges of the same color} \\ & \text{of\ } G \text{\ and\ } H \text{\ respectively} \\ 0 &\text{otherwise.} \endcases$$
Then $\QQ_{G,m}(B)$ is the number of maps $V \longrightarrow \{1, \ldots, k\}$ such that for every edge $\{u, v\}$ of $G$, the pair $\{\phi(u), \phi(v)\}$ spans an edge of $H$ of the same color and precisely $\mu_i$ vertices of $V$ are mapped into the vertex $i$ of $H$.

Clearly, (1.1.1) is a specialization of (1.5.1) since $\PP_{G,m}(A)=\QQ_{G,m}(B)$ provided 
$b^{uv}_{ij}=a_{ij}$ for all $\{u, v\} \in E$ and all $1 \leq i, j \leq k$. The advantage of working with $\QQ_{G,m}(B)$ instead of 
$\PP_{G,m}(A)$ is that $\QQ_{G,m}$ is a multi-affine polynomial, that is, the degree of each variable in $\QQ_{G,m}$ is 1.

If in (1.5.1) we consider the sum over all $k^{|V|}$ maps 
$\phi: V \longrightarrow \{1, \ldots, k\}$ we obtain the {\it partition function of edge-colored graph homomorphisms}
$$\QQ_G(B) =\sum_{\phi:\ V \rightarrow \{1, \ldots, k\}} \prod_{\{u, v\} \in E} b_{\phi(u) \phi(v)}^{uv}  \tag1.5.2$$
introduced in 
\cite{BS15}, cf. also \cite{AM98}.
\bigskip
$\bullet$ The main result of our paper is a deterministic algorithm, which, given a graph $G=(V, E)$, a vector 
$m=\left(\mu_1, \ldots, \mu_k\right)$ of multiplicities, a $|E| \times {k(k+1) \over 2}$ matrix $B=\left(b_{ij}^{uv}\right)$ such that 
$$\left| b^{uv}_{ij} - 1 \right| \ \leq \ {\gamma \over \Delta(G)} \quad \text{for all} \quad \{u, v\} \in E \quad \text{and all}
\quad 1 \leq  i, j \leq k  \tag1.5.3$$
and a real $\epsilon >0$, computes $\QQ_{G, m}(B)$ within relative error $\epsilon$ in $\left(|E|k\right)^{O(\ln |E| -\ln \epsilon)}$ time. Here $\gamma >0$ is an absolute constant, we can choose $\gamma=0.1$.
\bigskip
In \cite{BS15}, we construct a deterministic algorithm, which, given a graph $G=(V, E)$ and a matrix $B$ satisfying (1.5.3) with a better constant $\gamma=0.34$ (we can choose $\gamma=0.45$ 
if $\Delta(G) \geq 3$ and $\gamma =0.54$ for all sufficiently large $\Delta(G)$) computes the partition function (1.5.2)
within relative error $\epsilon$ in $\left(|E|k\right)^{O(\ln |E| -\ln \epsilon)}$ time. Although the methods of \cite{BS15} and this paper are similar, it appears that neither result follows from the other.

\subhead (1.6) Partition functions in combinatorics \endsubhead Partition functions are successfully used to count deterministically various combinatorial structures such as independent sets \cite{BG08} and matchings \cite{B+07} in graphs. 
The approach of \cite{BG08} and \cite{B+07} is based on the ``correlation decay" idea motivated by statistical physics.
Our approach is different (and the partition functions we compute are also different) but one can argue that our method is also inspired by statistical physics. Roughly, we use that the logarithm of the partition function is well-approximated by a low degree Taylor polynomial in the region that is sufficiently far away from the phase transition. Phase transitions are associated with {\it complex} zeros of partition functions \cite{LY52}, see also \cite{SS05} for connections to combinatorics, and the main effort of our method is in isolating the complex zeros of $\QQ_{G,m}$.

While our estimate of $\gamma=0.1$ in (1.1.3) and (1.5.3) is unlikely to be close to optimal, the tempting conjecture that $\PP_{G,m}(A)$ can be efficiently approximated as long as $\left| a_{ij}-1\right| < \gamma$ for any $\gamma < 1$, fixed in advance, is unlikely to be true even when $k=2$ and $\Delta(G)=3$. It is argued in \cite{BS15} that computing the partition function (1.1.2)  in quasi-polynomial time in such a large domain would have led to a quasi-polynomial algorithm in an NP-hard problem. The argument of \cite{BS15} transfers almost verbatim to the partition function (1.1.1). Hence it appears to be an interesting problem to find the best possible values of $\gamma$ in (1.1.3) and (1.5.3).   

\head 2. The algorithm \endhead

In this section, we describe the algorithm for computing (1.5.1). 
\subhead (2.1) The algorithm \endsubhead Let $J$ be an $|E| \times {k(k+1) \over 2}$ matrix filled with 1s. Given a
$|E| \times {k(k+1) \over 2}$ matrix $B=\left(b^{uv}_{ij}\right)$, where $\{u, v\} \in E$ and $1 \leq i,j \leq k$, we consider the univariate function 
$$f(t)=\ln \QQ_{G,m}\bigl(J +t(B-J)\bigr),$$
so that 
$$f(0)=\ln \QQ_{G,m}(J)=\ln {|V|! \over \mu_1! \cdots \mu_k!} \quad \text{and} \quad f(1)= \ln \QQ_{G,m}(B).$$
Hence our goal is to approximate $f(1)$ and we do it by using the Taylor polynomial approximation of $f$ at $t=0$:
$$f(1) \approx f(0) + \sum_{j=1}^n {1 \over j!} {d^j \over d t^j} f(t) \Big|_{t=0}. \tag2.1.1$$
We claim that the right hand side can be computed in $(|E|k)^{O(n)}$ time. Indeed, let 
$$\aligned g(t) =&\QQ_{G,m}\bigl(J + t(B-J)\bigr)\\=&\sum \Sb \phi:\ V \rightarrow \{1, \ldots, k\} \\ \left| \phi^{-1}(i)\right|=\mu_i \ \text{for} \ i=1, \ldots, k \endSb \prod_{\{u, v\} \in E} \left(1+t\left(b_{\phi(u)\phi(v)}^{uv}-1 \right)\right), \endaligned  \tag2.1.2$$
so that $f(t)=\ln g(t)$,
$$f'(t)={g'(t) \over g(t)} \quad \text{and} \quad g'(t) =g(t) f'(t).$$
Therefore,
$${d^j \over dt^j} g(t) \Big|_{t=0} = \sum_{i=0}^{j-1} {j-1 \choose i} \left({d^i \over dt^i} g(t) \Big|_{t=0}\right) \left({d^{j-i} \over dt^{j-i}} f(t) \Big|_{t=0}\right), \tag2.1.3$$
where we agree that the $0$-th derivative of $g$ is $g$. We note that 
$$g(0)={|V|! \over \mu_1! \cdots \mu_k!}. \tag2.1.4$$
If we compute the values of 
$${d^j \over dt^j} g(t) \Big|_{t=0} \quad \text{for} \quad j=1, \ldots, n, \tag2.1.5$$
then we can compute 
$${d^i \over dt^i} f(t) \Big|_{t=0} \quad \text{for} \quad i=1, \ldots, n, $$
from the triangular system (2.1.3) of linear equations with the coefficients (2.1.4)
on the diagonal. Hence our goal is to compute (2.1.5). 

Using (2.1.2), we obtain
$$\split &{d^j \over dt^j} g(t) \Big|_{t=0}  \\ \quad &= \sum \Sb \phi:\ V \rightarrow \{1, \ldots, k\} \\ \left| \phi^{-1}(i)\right|=\mu_i  \ \text{for} \ i=1, \ldots, k 
\endSb \sum \Sb I=\bigl( \{u_1, v_1\}, \\ \cdots \\  \{u_j, v_j\} \bigr) \endSb 
\left( b^{u_1 v_1}_{\phi(u_1) \phi(v_1)}-1 \right) \cdots \left(b^{u_j v_j}_{\phi(u_j) \phi(v_j)} -1 \right), \endsplit$$
where the inner sum is taken over all ordered collections $I$ of $j$ distinct edges $\{u_1, v_1\}$, $\ldots$, $\{u_j, v_j \}$ of $G$.
For such a collection $I$, let 
$$S(I) = \{u_1, v_1\} \cup \ldots \cup \{u_j, v_j\}$$
be the set of all distinct endpoints of the edges. Then we can write 
$$ \split {d^j \over dt^j} g(t)\Big|_{t=0} = &\sum \Sb I=\bigl( \{u_1, v_1\}, \\ \cdots \\  \{u_j, v_j\} \bigr) \endSb 
\sum \Sb \phi:\ S(I) \rightarrow \{1, \ldots, k\} \\ \left| \phi^{-1}(i) \right| \leq  \mu_i \text{\ for\ } i=1, \ldots, k \endSb 
{\bigl(|V| -|S(I)|)! \over \left(\mu_1 - \left| \phi^{-1}(1)\right| \right)! \cdots \left(\mu_k - \left| \phi^{-1}(k)\right| \right)!} \\
&\quad \times \left( b^{u_1 v_1}_{\phi(u_1) \phi(v_1)}-1 \right) \cdots \left(b^{u_j v_j}_{\phi(u_j) \phi(v_j)} -1 \right).\endsplit$$
In words: we enumerate at most $|E|^j$ ordered collections $I=\bigl( \left\{u_1, v_1\right\}, \ldots$, \break $\left\{u_j, v_j \right\} \bigr)$ of $j$ distinct edges of $G$, for each such a collection, we enumerate at most $k^{2j}$ maps $\phi$ defined on the set of the endpoints of the edges from $I$ into the set $\{1, \ldots, k\}$, 
multiply the term  $\left( b^{u_1 v_1}_{\phi(u_1) \phi(v_1)}-1 \right) \cdots \left(b^{u_j v_j}_{\phi(u_j) \phi(v_j)} -1 \right)$ by the number of ways to extend the map $\phi$ to the whole set $V$ of vertices so that the multiplicity of $i$ is $\mu_i$ 
and add the results over all choices of $I$ and $\phi$. Since $j \leq n$, the complexity of computing (2.1.5) is indeed 
$(|E| k)^{O(n)}$ as claimed.

The quality of the approximation (2.1.1) depends on the location of {\it complex} zeros of $\QQ_{G,m}$.
\proclaim{(2.2) Lemma} Suppose that there is a real $\beta >1$ such that 
$$\QQ_{G,m}\bigl(J + z(B-J)\bigr) \ne 0 \quad \text{for all} \quad z \in {\Bbb C} \quad \text{satisfying} \quad |z| \leq \beta.$$
Then the right hand side of (2.1.1) approximates $f(1)$ within an additive error of 
$${|E| \over (n+1) \beta^n (\beta -1)}.$$
\endproclaim 
\demo{Proof} The function $g(t)$ defined by (2.1.2) is a polynomial in $t$ of degree at most $d \leq |E|$ and $g(0)\ne 0$
by (2.1.4),
so we can factor
$$g(z)=g(0) \prod_{i=1}^d \left(1-{z \over \alpha_i}\right),$$
where $\alpha_1, \ldots, \alpha_d \in {\Bbb C}$ are the roots of $g(z)$. In addition,
$$\left| \alpha_i \right| \ \geq \ \beta \ > \ 1 \quad \text{for} \quad i=1, \ldots, d.$$
Thus 
$$f(z) =\ln g(z) =\ln g(0) + \sum_{i=1}^d \ln \left(1-{z \over \alpha_i}\right) \quad \text{for all} \quad |z| \leq 1, \tag2.2.1$$
where we choose the branch of $\ln g(z)$ for which $\ln g(0)$ is real. Using the standard Taylor series expansion, we obtain
$$\ln \left(1 -{1 \over \alpha_i}\right) =-\sum_{j=1}^n {1 \over j} \left({1 \over \alpha_i}\right)^j +\xi_n,$$
where 
$$\left| \xi_n\right| =\left| \sum_{j=n+1}^{+\infty} {1 \over j} \left({1 \over \alpha_i}\right)^j \right| \ \leq \ {1 \over (n+1) \beta^n (\beta-1)}.$$
Hence from (2.2.1) we obtain
$$f(1)=f(0) +\sum_{j=1}^n \left(-{1 \over j} \sum_{i=1}^d \left({1 \over \alpha_i}\right)^j \right) +\eta_n,$$
where 
$$\left| \eta_n \right| \ \leq \ {|E| \over (n+1) \beta^n (\beta-1)}.$$
It remains to notice that 
$$-{1\over j} \sum_{i=1}^d \left({1 \over \alpha_i}\right)^j ={1 \over j!} {d^j \over dt^j} f(t) \Big|_{t=0}.$$
{\hfill \hfill \hfill} \qed
\enddemo

For a fixed $\beta >1$, to achieve an additive error of $0 < \epsilon <1$, we can choose $n=O\left(\ln |E| -\ln \epsilon\right)$, in which case the algorithm of Section 2.1 computes $\QQ_{G,m}(B)$ within a relative error $\epsilon$ in $(|E| k)^{O(\ln |E|-\ln \epsilon)}$ time. Hence it remains to identify matrices $B$ for which the number $\beta >1$ of Lemma 2.2 exists.

We prove the following result.

\proclaim{(2.3) Theorem} There exists an absolute constant $\alpha >0$ (one can choose $\alpha=0.107$) such that for any undirected graph $G=(V, E)$, for any vector \newline $m=\left(\mu_1, \ldots, \mu_k\right)$ of positive integers such that $\mu_1 +\ldots +\mu_k =|V|$ and for any complex 
$|E| \times {k(k+1)\over 2}$ matrix $Z=\left(z^{uv}_{ij}\right)$ satisfying 
$$\left| 1- z^{uv}_{ij} \right| \ \leq \ {\alpha \over \Delta(G)} \quad \text{for all} \quad \{u, v\} \in E \quad \text{and} \quad 
1 \leq i, j \leq k,$$
one has 
$$\QQ_{G,m}(Z) \ne 0.$$
\endproclaim

Theorem 2.3 implies that if $B$ satisfies (1.5.3) with $\gamma=0.100$, then we can choose 
$$\beta= {\alpha \over \gamma}={0.107 \over 0.1} =1.07 > 1$$ 
in Lemma 2.2 and hence we obtain an algorithm which computes $\QQ_{G,m}(B)$ within relative error $\epsilon$ in 
$(|E| k)^{O(\ln |E| -\ln \epsilon)}$ time.

In the rest of the paper, we prove Theorem 2.3.

\head 3. Recurrence relations \endhead 

We consider the polynomials $\QQ_{G,m}(Z)$ within a larger family of polynomials.

\subhead (3.1) Definitions \endsubhead
Let us fix the graph $G$ and the integer vector $m=\left(\mu_1, \ldots, \mu_k\right)$.

We say that a sequence $W=\left(v_1, \ldots, v_n\right)$ of vertices of $G$ is {\it admissible} if the vertices $v_1, \ldots, v_n$ are distinct.

Let $I=\left(i_1, \ldots, i_n\right)$ be a sequence of (not necessarily distinct) indices $i_j \in \{1, \ldots, k\}$ for $j=1, \ldots, n$. For $i \in \{1, \ldots, k\}$, we define the {\it multiplicity $\nu_i(I)$ of $i$ in $I$} by 
$$\nu_i(I)=\left|\left\{j:\ i_j=i \right\}\right|,$$ the number of times that $i$ occurs in $I$. 

A sequence $I=\left(i_1, \ldots, i_n\right)$ of (not necessarily distinct) indices $i_j \in \{1, \ldots, k\}$ is called {\it admissible} 
provided
$$\nu_i(I) \ \leq \ \mu_i \quad \text{for} \quad i=1,\ldots, k.$$
For an admissible sequence $W=\left(v_1, \ldots, v_n\right)$  of vertices and an admissible sequence $I=\left(i_1, \ldots, i_n\right)$ of indices such that $|W|=|I|$, we define
$$\QQ^W_I(Z)=\sum \Sb \phi: V \rightarrow \{1, \ldots, k\} \\ \left| \phi^{-1}(i)\right|=\mu_i \ \text{for}\ i=1, \ldots, k 
\\ \phi\left(v_j\right)=i_j \ \text{for} \quad j=1, \ldots, n \endSb 
\prod_{\{u, v\} \in E} z^{uv}_{\phi(u)\phi(v)}.$$
In other words, to define $\QQ^W_I(Z)$, we restrict the sum (1.5.1) defining $\QQ_{G,m}(Z)$ to maps $\phi: V \longrightarrow \{1, \ldots, k\}$ that map prescribed vertices $v_1, \ldots, v_n$ of the graph to the prescribed indices $i_1, \ldots, i_n$. Note that 
if $W=\left(v_1, \ldots, v_n\right)$ and $I=\left(i_1, \ldots, i_n\right)$ are admissible sequences, then there is a map $\phi: V \longrightarrow \{1, \ldots, k\}$ which maps precisely $\mu_i$ distinct vertices of $V$ onto $i$ for $i=1, \ldots, k$ and such that $\phi(v_j)=i_j$ for $j=1, \ldots, n$.

We suppress the graph $G$ and the vector of multiplicities $m$ in the notation for $\QQ^W_I(Z)$ and note that when $W$ and $I$ are both empty, then $\QQ^W_I(Z)=\QQ_{G,m}(Z)$.

For a sequence $W$ of vertices and a vertex $v$ of the graph, we denote by $(W, v)$ the sequence $W$ appended by $v$.
For a sequence $I$ of indices and another index $i \in \{1, \ldots, k\}$, we denote by $(I, i)$ the sequence $I$ appended by $i$. We denote similarly sequences appended by several vertices or indices.

\subhead (3.2) Recurrence relations \endsubhead We will use the following two recurrence relations.

First, if $W$ and $I$ are admissible sequences such that $|W|=|I|$ and the sequence 
$(W, v)$ is also admissible (that is, $v$ is distinct from the vertices in $W$), then 
$$\QQ^W_I(Z) =\sum \Sb i \in \{1, \ldots, k\}: \\ (I, i) \ \text{is admissible} \endSb \QQ^{(W, v)}_{(I, i)}(Z). \tag3.2.1$$

Second, if $W$ and $I$ are admissible sequences such that $|W|=|I|$ and 
if $(I, i)$ is an admissible sequence of indices then 
$$\QQ^W_I(Z)= {1 \over \mu_i - \nu_i(I)}\sum \Sb v \in V: \\ (W, v) \ \text{is admissible} \endSb Q^{(W, v)}_{(I, i)}(Z) \tag 3.2.2$$
(recall that $\nu_i(I)$ is the multiplicity of $i$ in $I$ so that $\mu_i -\nu_i(I) > 0$ if $(I, i)$ is admissible).

\head 4. Angles in the complex plane \endhead

In what follows, we measure angles between non-zero complex numbers, considered as vectors in ${\Bbb R}^2$ identified with ${\Bbb C}$.

We start with two geometric calculations.

\proclaim{(4.1) Lemma} Let $z_1, \ldots, z_n \in {\Bbb C}$ be non-zero numbers such that the angle between any two numbers $z_i$ and $z_j$ does not exceed $\theta$ for some $0 \leq \theta < 2\pi/3$. Then, for $z=z_1+ \ldots + z_n$, we have 
$$|z| \ \geq \ \left(\cos {\theta \over 2}\right) \sum_{i=1}^n |z_i|.$$
\endproclaim
\demo{Proof} We observe that $0$ is not in the convex hull of any three vectors $z_i, z_j, z_k$, since otherwise the angle between some two of those three vectors would have been at least $2 \pi/3$. By the Carath\'eodory Theorem, $0$ is not in the convex hull of $z_1, \ldots, z_n$. Therefore, the vectors are enclosed in an angle of at most $\theta$ with vertex at $0$. Let us project each vector $z_i$ orthogonally onto the bisector of the angle. The length of the projection is at least $|z_i| \cos(\theta/2)$ and hence the length of the projection of $z_1 + \ldots + z_n$ is at least $(|z_1| + \ldots + |z_n|)\cos(\theta/2)$. Since the length of a vector is at least the length of its orthogonal projection, the proof follows.
{\hfill \hfill \hfill} \qed
\enddemo

Lemma 4.1 was suggested by Boris Bukh \cite{Bu15}. It replaces a weaker bound of 
$\sqrt{\cos \theta}\left(|z_1| + \ldots +|z_n|\right)$ of an earlier version of the paper.

\proclaim{(4.2) Lemma} Let $a_1, \ldots, a_n$ and $b_1, \ldots, b_n$ be complex numbers such that all $a_1, \ldots, a_n$ are non-zero. 
Let 
$$a=\sum_{i=1}^n a_i \quad \text{and} \quad b= \sum_{i=1}^n b_i.$$
Suppose that for some real $1> \tau > \epsilon >0$ we have 
$$\split &\left| {b_i \over a_i} -1 \right| \ \leq \ \epsilon \quad \text{for} \quad i=1, \ldots, n \qquad \text{and} \\
&|a| \ \geq \ \tau \sum_{i=1}^n \left| a_i\right|. \endsplit$$
Then $a \ne 0$, $b\ne 0$ and the angle between $a$ and $b$ does not exceed 
$$\arcsin {\epsilon \over \tau}.$$
\endproclaim
\demo{Proof}  Clearly, $a \ne 0$. We can write
$$b_i=\left(1+\epsilon_i\right)a_i \quad \text{where} \quad \left| \epsilon_i \right| \ \leq \ \epsilon \quad \text{for} \quad i=1, \ldots, n.$$
Hence 
$$b =\sum_{i=1}^n \left(1+\epsilon_i\right) a_i = a + \sum_{i=1}^n \epsilon_i a_i \quad \text{and} \quad 
\left| \sum_{i=1}^n \epsilon_i a_i \right| \ \leq \ \epsilon \sum_{i=1}^n \left| a_i \right| \ \leq \ {\epsilon \over  \tau} |a|,$$
from which $b \ne 0$ and 
$$\left|{b \over a} -1 \right| \ \leq \ {\epsilon \over \tau}.$$
Therefore, the argument of the complex number $b/a$ lies in the interval 
$$\left[-\arcsin {\epsilon\over \tau},\ \arcsin {\epsilon\over \tau} \right]$$ and the proof follows.
{\hfill \hfill \hfill} \qed
\enddemo

The main result of this section concerns the angles between various numbers $\QQ^W_I(Z)$ introduced in Section 3.

\proclaim{(4.3) Proposition} Let us us fix an admissible sequence $W$ of vertices, an admissible sequence $I$ of indices such that $0 \leq |W|=|I| \leq |V|-2$, a complex $|E| \times {k(k+1) \over 2}$ matrix $Z$, a real $\epsilon >0$ and a real
$0 \leq \theta < 2\pi/3$ such that $\epsilon < \cos( \theta/2)$. Suppose that for any two vertices $u, v \in V$ and any two
$i, j \in \{1, \ldots, k\}$ such that the sequences $(W, u, v)$ and $(I, i, j)$ are admissible, we have 
$\QQ^{(W, u, v)}_{(I, i, j)}(Z) \ne 0$,  $\QQ^{(W, u, v)}_{(I, j, i)}(Z) \ne 0$ and
$$\left| {\QQ^{(W, u, v)}_{(I, i, j)}(Z) \over \QQ^{(W, u, v)}_{(I, j, i)}(Z)} -1 \right| \ \leq \ \epsilon.$$
\roster
 \item Suppose that for any two vertices $u, v \in V$ and any three 
$i, j_1, j_2 \in \{1, \ldots, k\}$ such that the sequences $(W, u, v)$ and $(I, i, j_1)$ and $(I, i, j_2)$ are admissible, the angle between two complex numbers $\QQ^{(W, u, v)}_{(I, i, j_1)}(Z)$ and $\QQ^{(W, u, v)}_{(I, i, j_2)}(Z)$ does not exceed $\theta$.
Then for any two vertices $u, v \in V$ and for any $i \in \{1, \ldots, k\}$ such that the sequences $(W, u)$, $(W, v)$ and $(I, i)$ are admissible, we have $\QQ^{(W, u)}_{(I, i)}(Z) \ne 0$, $\QQ^{(W, v)}_{(I, i)}(Z) \ne 0$ and the angle between two complex numbers $\QQ^{(W, u)}_{(I, i)}(Z)$ and $\QQ^{(W, v)}_{(I, i)}(Z)$ does not exceed
$$\arcsin {\epsilon \over \cos(\theta/2)}.$$
\item Suppose that for any three $u, v_1, v_2 \in V$ and any two indices  $i, j \in \{1, \ldots, k\}$ such that the sequences
$(W, u, v_1)$, $(W, u, v_2)$ and $(I, i, j)$ are admissible, the angle between two complex numbers $\QQ^{(W, u, v_1)}_{(I, i, j)}(Z)$ and $\QQ^{(W, u, v_2)}_{(I, i, j)}(Z)$ does not exceed $\theta$. Then for any $u \in V$ and any $i, j \in \{1, \ldots, k\}$ such that the sequences $(W, u)$, $(I, i)$ and $(I, j)$ are admissible, we have $\QQ^{(W, u)}_{(I, i)}(Z) \ne 0$, $\QQ^{(W, u)}_{(I, j)}(Z) \ne 0$ and the angle between two complex numbers $\QQ^{(W, u)}_{(I, i)}(Z)$ and $\QQ^{(W, u)}_{(I, j)}(Z)$ does not exceed 
$$\arcsin {\epsilon \over \cos(\theta/2)}.$$
\endroster
\endproclaim

\demo{Proof} To prove Part (1), we note that by (3.2.1), we have
$$\split \QQ^{(W, u)}_{(I, i)}(Z)=&\sum \Sb j \in \{1, \ldots, k\}: \\ (I, i, j) \ \text{is admissible} \endSb \QQ^{(W, u, v)}_{(I, i, j)}(Z)
\quad \text{and} \quad \\
\QQ^{(W, v)}_{(I, i)}(Z)=&\sum \Sb j \in \{1, \ldots, k\}: \\ (I, j, i) \ \text{is admissible} \endSb \QQ^{(W, v, u)}_{(I, i, j)}(Z).\endsplit$$
Since the angle between any two complex numbers $\QQ^{(W, u, v)}_{(I, i, j_1)}(Z)$ and $\QQ^{(W, u, v)}_{(I, i, j_2)}(Z)$ does not exceed $\theta$, by Lemma 4.1 we have 
$$\left| \QQ^{(W, u)}_{(I, i)}(Z) \right| \ \geq \ \left(\cos{\theta \over 2}\right) \sum \Sb j \in \{1, \ldots, k\}: \\ (I, i, j) \ \text{is admissible} \endSb 
\left| \QQ^{(W, u, v)}_{(I, i, j)}(Z) \right|.$$
The proof of Part (1) follows by Lemma 4.2 applied to the numbers
$$a_j=\QQ^{(W, u, v)}_{(I, i, j)}(Z), \quad b_j =\QQ^{(W, v, u)}_{(I, i, j)}(Z), \quad a = \QQ^{(W, u)}_{(I, i)}(Z), \quad 
 b=\QQ^{(W, v)}_{(I, i)}(Z)$$ and 
 $\tau=\cos(\theta/2)$.

To prove Part (2), let us fix two indices $i \ne j$. Then the sequence $(I, i, j)$ is admissible provided $(I, i)$ and $(I, j)$ are both admissible. Applying (3.2.2), we obtain
$$\split \QQ^{(W, u)}_{(I, i)}(Z)=&{1 \over \mu_j -\nu_j(I,i)} \sum \Sb v \in V: \\ (W, u,v) \ \text{is admissible} \endSb 
\QQ^{(W, u, v)}_{(I, i, j)}(Z) \\ &\qquad \qquad \text{and} \\ 
\QQ^{(W, u)}_{(I, j)}(Z)=&{1 \over \mu_i -\nu_i(I, j)} \sum \Sb v \in V: \\ (W,u,v) \ \text{is admissible} \endSb 
\QQ^{(W, u, v)}_{(I, j, i)}(Z).\endsplit \tag4.3.1
$$
For a vertex $v \in V$ such that the sequence $(W, u, v)$ is admissible, let us denote
$$a_v =\QQ^{(W, u, v)}_{(I, i, j)}(Z) \quad \text{and} \quad b_v=\QQ^{(W, u, v)}_{(I, j, i)}(Z)$$
and let
$$a=\sum \Sb v \in V: \\ (W,u,v) \text{\ is admissible} \endSb a_v \quad \text{and} \quad b=\sum \Sb v \in V:\\ (W,u,v)\ \text{is admissible} \endSb b_v.$$
Since the angle between any two numbers $a_{v_1}$ and $a_{v_2}$ does not exceed $\theta$, by Lemma 4.1 we have 
$$|a| \ \geq \ \left(\cos{\theta \over 2}\right)\sum \Sb v \in V: \\ (W,u,v) \text{\ is admissible} \endSb \left| a_v \right|,$$
and by Lemma 4.2, the angle between $a$ and $b$ does not exceed 
$$\arcsin {\epsilon \over \cos(\theta/2)}.$$ Since by (4.3.1) we have 
$$ \QQ^{(W, u)}_{(I, i)}(Z) ={1 \over \mu_j -\nu_j(I,i)} a \quad \text{and} \quad \QQ^{(W, u)}_{(I, j)}(Z)={1 \over \mu_i -\nu_i(I, j)} b,$$
where $\mu_j -\nu_j(I, i)>0$ and $\mu_i -\nu_i(I, j)>0$, the result follows.

It remains to handle the case of $i=j$, in which case we just have to prove that $\QQ^{(W, u)}_{(I, i)}(Z) \ne 0$ whenever 
$(W, u)$ and $(I, i)$ are admissible sequences. Since $|W|=|I| \leq |V|-2$, there is an index $l \in \{1, \ldots, k\}$ for which the 
sequence $(I, i, l)$ is admissible. Using the first equation in (4.3.1) with $j$ replaced by $l$ throughout and Lemma 4.1, we conclude that $\QQ^{(W, u)}_{(I, i)}(Z) \ne 0$.
{\hfill \hfill \hfill} \qed
\enddemo

\remark{(4.4) Remark} The intuition behind Proposition 4.3 is as follows. We are interested in how the value of $\QQ^W_I(Z)$ may change if we change one vertex in the sequence $W$ or one index in the sequence $I$, while keeping the sequences admissible. 
Suppose that the value of $\QQ^W_I(Z)$ does not change much if we just permute two indices in the sequence $I$. Part (1) asserts that if for sequences $W$ and $I$ of some length $|W|=|I|$, the complex number $\QQ^W_I(Z)$ rotates by a small angle when one index in $I$ is changed, then for {\it shorter} sequences $W'$ and $I'$ such that $|W'|=|W|-1$ and 
$|I'|=|I|-1$, the complex number $\QQ^{W'}_{I'}(Z)$ rotates by a small angle when one vertex in $W'$ is changed. Part (2) asserts that if for sequences $W$ and $I$ of some length $|W|=|I|$, the complex number $\QQ^W_I(Z)$ rotates by a small angle when one vertex in $W$ is changed, then for {\it shorter} sequences $W'$ and $I'$ such that $|W'|=|W|-1$ and 
$|I'|=|I|-1$, the complex number $\QQ^{W'}_{I'}(Z)$ rotates by a small angle when one index in $I'$ is changed.

For the subsequent proof, we would like this condition of $\QQ^W_I(Z)$ being rotated by a small angle if one vertex in $W$ is changed or one index in $I$ is changed to propagate for shorter and shorter sequences $W$ and $I$. For that, we will need to find the ``fixed point" of the conditions of Proposition 4.3, that is, a solution $0 \leq \theta < 2\pi/3$ to the equation 
$$\theta = \arcsin {\epsilon \over \cos(\theta/2)}.$$
It is not hard to see that for all sufficiently small $\epsilon >0$ there is such a solution. To make the constant $\alpha$ in Theorem 2.3 as large as possible, we would like to choose $\epsilon$ as large as possible. Numerical computations show that we can choose $\epsilon=0.76$, in which case $\theta \approx 1.101463960$.
\endremark

\head 5. Derivatives \endhead

The first goal of this section is to relate how the value of $\QQ^W_I(Z)$ changes when two indices in $I$ are permuted, with partial derivatives of $\QQ^W_I(Z)$.

\definition{(5.1) Definition} For a $0 < \delta <1$, we define the polydisc $\UU(\delta) \subset {\Bbb C}^{k(k+1)|E|/2}$ by
$$\UU(\delta)=\left\{ Z=\left(z^{uv}_{ij}\right): \quad \left|1-z^{uv}_{ij}\right| \leq \delta \quad \text{for all} \quad 
\{u, v\} \in E \quad \text{and} \quad 1 \leq i, j \leq k \right\}.$$
\enddefinition

Thus $\UU(\delta)$ is the closed polydisc of radii $\delta$ centered at the matrix $J$ of all 1's. We will be interested in the situations when $\QQ^W_I(Z) \ne 0$ for all $Z \in \UU(\delta)$, in which case we can choose a branch of $\ln \QQ^W_I(Z)$ for $Z \in \UU(\delta)$ in such a way that $\ln \QQ^W_I(J)$ is a real number.

\proclaim{(5.2) Lemma} Let us fix an integer $2 \leq r \leq |V|$ and let $\tau >0$ and $0 < \delta < 1$ be real. Suppose that 
for any admissible sequences $W$ and $I$ such that $|W|=|I|=r$ and for any $Z \in \UU(\delta)$ we have 
$\QQ^W_I(Z) \ne 0$ and the following condition is satisfied: if $W=(W', v)$ and $I=(I', i)$ then 
$$\left|\QQ^W_I(Z) \right| \ \geq \ {\tau \over \Delta(G)} \sum \Sb w:\ \{w, v\} \in E \\ l:\ 1 \leq l \leq k \endSb 
\left|z^{vw}_{il} \right| \left| {\partial \over \partial z^{vw}_{il}} \QQ^W_I(Z) \right|.$$

Then for any admissible $W$ and $I$ such that $|W|=|I|=r$ and for any $Z \in \UU(\delta)$ the following condition is satisfied:
if $W=(W',u,v)$ and $I=(I', j, i)$, then for
$$\xi={4 \delta \Delta(G) \over (1-\delta) \tau}$$ we have 
$$\left| {\QQ^{(W', u, v)}_{(I', j, i)}(Z) \over \QQ^{(W', u, v)}_{(I',i,j)}(Z)} -1 \right| \ \leq \ e^{\xi}-1.$$
\endproclaim
\demo{Proof} Let us choose admissible $W$ and $I$ such that $|W|=|I|=r$ and suppose that 
$W=(W', u, v)$ and $I=(I', j, i)$. If $i=j$ the result is trivial, so we assume that $i \ne j$. Since 
$${\partial \over \partial z^{vw}_{il}} \ln \QQ^W_I(Z)= \left({\partial \over \partial z^{vw}_{il}} \QQ^W_I(Z)\right)/\QQ^W_I(Z)$$
and 
$$\left|z^{vw}_{il} \right| \ \geq \ 1-\delta \quad \text{for all} \quad Z \in \UU(\delta),$$
by the conditions of the lemma, we obtain 
$$\split &\sum \Sb w:\ \{w, v\} \in E \\ l:\ 1 \leq l \leq k \endSb \left| {\partial \over \partial z_{il}^{vw}} \ln \QQ^W_I(Z) \right| \ \leq \ 
{\Delta(G) \over (1-\delta) \tau} \quad \text{and} \quad \\ 
&\sum \Sb w:\ \{w, u\} \in E \\ l:\ 1 \leq l \leq k \endSb \left| {\partial \over \partial z_{jl}^{uw}} \ln \QQ^W_I(Z) \right| \ \leq \ 
{\Delta(G) \over (1-\delta) \tau}.  \endsplit \tag5.2.1$$
Given a matrix $A \in \UU(\delta)$, we define a matrix $B \in \UU(\delta)$ as follows:
$$\split &b^{uw}_{jl} = a^{uw}_{il} \quad \text{for all} \quad w\ne v \quad \text{such that} \quad \{u, w\} \in E \quad \text{and all} \quad l=1, \ldots, k \\
&b^{vw}_{il}=a^{vw}_{jl} \quad \text{for all} \quad w\ne u \quad \text{such that} \quad \{v, w\} \in E \quad \text{and all} \quad 
l=1, \ldots, k, \endsplit$$
whereas all other entries of $B$ are equal to the corresponding entries of $A$.
Since swapping the values of $\phi: V \longrightarrow \{1, \ldots, k\}$ on two vertices $u$ and $v$ does not change the multiplicities $\left| \phi^{-1}(l) \right|$ for $l=1, \ldots, k$, we have
$$\QQ^{(W', u, v)}_{(I',j,i)}(B)=\QQ^{(W', u, v)}_{(I',i,j)}(A)$$
and using (5.2.1), we obtain
$$\split &\left| \ln \QQ^{(W', u, v)}_{(I', j, i)}(A) - \ln \QQ^{(W', u, v)}_{(I', i, j)}(A)\right| =
\left| \ln \QQ^{(W', u, v)}_{(I', j, i)}(A) - \ln \QQ^{(W', u, v)}_{(I', j,  i)}(B) \right| \\ &\quad \leq \ 
\left(\sup_{Z \in \UU(\delta)} 
 \sum \Sb w:\ \{u, w\} \in E \\ l:\ 1 \leq l \leq k \endSb \left| {\partial \over \partial z^{uw}_{jl}} \ln \QQ^W_I(Z) \right|
+ \sum \Sb w:\ \{v, w\} \in E \\ l:\ 1 \leq l \leq k \endSb \left| {\partial \over \partial z^{vw}_{il}} \ln \QQ^W_I(Z) \right| \right) 
\\ &\qquad \qquad \times \left( \max \Sb w \in V \\ 1 \leq l \leq k \endSb \left| a^{uw}_{jl}-b^{uw}_{jl} \right|, \ \left| a^{vw}_{il}-b^{vw}_{il} \right| \right) \\
&\leq \ \left({\Delta(G) \over (1-\delta) \tau} + {\Delta(G) \over (1-\delta) \tau}\right) \times (2\delta)= {4 \delta \Delta(G) \over (1-\delta) \tau} =\xi.
\endsplit $$
Let us denote 
$$\zeta={\QQ^{(W', u, v)}_{(I', j, i)}(Z) \over \QQ^{(W', u, v)}_{(I',i,j)}(Z)}.$$
Hence 
$$|\ln \zeta|  \ \leq \ \xi.$$
Since 
$$\left| e^z -1 \right| =\left| \sum_{m=1}^{+\infty} {z^m \over m!} \right| \ \leq \ \sum_{m=1}^{+\infty} {|z|^m \over m!} =e^{|z|} -1 \quad 
\text{for all} \quad z \in {\Bbb C},$$
we conclude that 
$$\left| \zeta - 1 \right| \ \leq \ \left| e^{\ln \zeta} -1 \right| \ \leq \ e^{|\ln \zeta|} -1 \ \leq \ e^{\xi} -1$$
as desired.
{\hfill \hfill \hfill} \qed
\enddemo

\remark{(5.3) Remark} We will be interested in the situation when 
$$\left| {\QQ^{(W', u, v)}_{(I', j, i)}(Z) \over \QQ^{(W', u, v)}_{(I',i,j)}(Z)} -1 \right| \ \leq \ 0.76,$$
cf. Proposition 4.3 and Remark 4.4. To ensure the estimate, it suffices to have $\xi \leq 0.565$ in Lemma 5.2.
\endremark

Our next goal is to relate the inequality for partial derivatives of Lemma 5.2 to another angle condition of Proposition 4.3. Namely, we show that if we can bound the angle by which the number $\QQ^W_I(Z)$ rotates when one index in $I$ is changed, we obtain the inequality of Lemma 5.2 for {\it shorter} sequences $W'$ and $I'$ such that $|W'|=|I'|=|W|-1=|I|-1$.

\proclaim{(5.4) Lemma} Let $0 \leq \theta < 2\pi/3$ be a real number, let $W$ be an admissible sequence of vertices and let $I$ be an admissible sequence of indices such that $1 \leq |I|=|W| \leq |V|-1$. Suppose that for any $Z \in \UU(\delta)$, for every $w$ such that 
$(W, w)$ is admissible and for every $1 \leq l, j \leq k$ such that $(I, l)$ and $(I, j)$ are admissible, we have 
$\QQ^{(W, w)}_{(I, l)}(Z) \ne 0$, $\QQ^{(W, w)}_{(I, j)}(Z) \ne 0$ and the angle between the complex numbers 
$\QQ^{(W, w)}_{(I, l)}(Z)$ and $\QQ^{(W, w)}_{(I, j)}(Z)$  does not exceed $\theta$.

Let $W=(W', v)$ and $I=(I', i)$.
Then 
$$\left| \QQ^W_I(Z)\right| \ \geq \ {\tau \over \Delta(G)} \sum \Sb w:\ \{w, v\} \in E \\ j:\ 1 \leq j \leq k \endSb 
\left|z^{vw}_{ij}\right| \left| {\partial \over \partial z^{vw}_{ij}} \QQ^W_I(Z) \right|,$$
where 
$$\tau=\cos {\theta \over 2}.$$
\endproclaim
\demo{Proof} Recall that $\QQ^W_I(Z)$ is a multi-affine function in $Z$.

Let us choose a vertex $w$ such that $\{v, w\} \in E$.
If $w$ is an element of the sequence $W'$, then
$$z^{vw}_{ij} {\partial \over \partial z^{vw}_{ij}}\QQ^W_I(Z) =\cases \QQ^W_I(Z) &\text{if $j$ is the element of $I'$ that corresponds to $w$} \\ 0 &\text{otherwise.}
\endcases$$
If $w$ is not an element of $W'$, then $(W, w)$ is an admissible sequence and 
$$z^{vw}_{ij} {\partial \over \partial z^{vw}_{ij}}\QQ^W_I(Z) =\cases 
\QQ^{(W, w)}_{(I, j)}(Z) &\text{if $(I, j)$ is admissible,} \\ 0 &\text{otherwise.} \endcases$$
Denoting by $d_0$ the number of vertices $w$ in $W'$ such that $\{w, v\} \in E$, we obtain
$$\split &\sum \Sb w:\ \{w, v\} \in E \\ j:\ 1 \leq j \leq k \endSb 
\left|z^{vw}_{ij}\right| \left| {\partial \over \partial z^{vw}_{ij}} \QQ^W_I(Z) \right| =
d_0 \left|\QQ^W_I(Z)\right| \\ &\qquad + \sum \Sb w, j: \\ w \text{\ not in $W'$}, \{w, v\} \in E \\ (I, j) \text{\ is admissible} \endSb 
\left| \QQ^{(W, w)}_{(I, j)}(Z)\right|. \endsplit \tag5.4.1$$
On the other hand, by (3.2.1) for all $w$ not in $W'$, we have 
$$\QQ^W_I(Z) =\sum \Sb j:\ 1 \leq j \leq k \\ (I, j) \text{\ is admissible} \endSb \QQ^{(W, w)}_{(I, j)}(Z)$$ and hence by Lemma 4.1, 
$$\left|\QQ^W_I(Z) \right| \ \geq \ \left(\cos{\theta \over 2}\right) \sum \Sb j:\ 1 \leq j \leq k \\ (I, j) \text{\ is admissible} \endSb 
\left| \QQ^{(W, w)}_{(I, j)}(Z) \right|. \tag5.4.2$$
Denoting by $d_1$ the number of vertices $w$ such that $\{v, w\} \in E$ and $w$ are not in $W'$, we deduce from (5.4.1) and (5.4.2)
that 
$$\split &{\cos(\theta/2) \over \Delta(G)} \sum \Sb w:\ \{w, v\} \in E \\ j:\ 1 \leq j \leq k \endSb 
\left|z^{vw}_{ij}\right| \left| {\partial \over \partial z^{vw}_{ij}} \QQ^W_I(Z) \right|\\ &\quad = {d_0 \cos(\theta/2) \over \Delta(G)} \left|\QQ^W_I(Z)\right|+{\cos(\theta/2) \over \Delta(G)} \sum \Sb w, j: \\ w \text{\ not in $W'$}, \{w, v\} \in E \\ (I, j) \text{\ is admissible} \endSb 
\left| \QQ^{(W, w)}_{(I, j)}(Z)\right| \\ &\quad \leq  {d_0 \cos(\theta/2) \over \Delta(G)} \left|\QQ^W_I(Z)\right| 
+{d_1 \over \Delta(G)} \left|\QQ^W_I(Z)\right| \leq \left|\QQ^W_I(Z)\right|,\endsplit$$
as desired.
{\hfill \hfill \hfill} \qed
\enddemo

\head 6. Proof of Theorem 2.3 \endhead

First, we define some constants. Let
$$\epsilon =0.76$$ 
and let $\theta$ be the solution of the equation 
$$\theta=\arcsin {\epsilon \over \cos( \theta/2)},$$
so,
$$ \theta \approx 1.101463960,$$
cf. Remark 4.4. Let
$$\tau=\cos{\theta \over 2} \approx 0.8521416971$$
and let
$$\alpha = {0.565 \tau \over 4+ 0.565 \tau} \approx  0.1074337498.$$
In particular, as long as $\Delta(G) \geq 1$, we have 
$$\xi={4 \alpha \over \left(1- {\alpha \over \Delta(G)}\right) \tau} \ \leq \ 0.565,$$
and hence
$$e^{\xi}-1  \ \leq \ \epsilon,$$
cf. Remark 5.3.
Finally, let 
$$\delta = {\alpha \over \Delta(G)}.$$
Recall that $\UU(\delta)$ is the polydisc of radii $\delta$ centered at the matrix of all 1's, cf. Definition 5.1.

We prove by descending induction for $r=|V|, |V|-1, \ldots, 2$ that the following statements (6.1.$r$)--(6.5.$r$) are satisfied for all
$Z \in \UU(\delta)$.
\bigskip
(6.1.$r$) Let $W$ be an admissible sequence of vertices and let $I$ be an admissible sequence of indices such that 
$|W|=|I|=r$. Then $\QQ^W_I(Z) \ne 0$.
\medskip 
(6.2.$r$) Let $W$ be an admissible sequence of vertices and let $I$ be an admissible sequence of indices such that 
$|W|=|I|=r$. Assuming that $W=(W', v)$ and $I=(I', i)$, we have
$$\left| Q^W_I(Z)\right| \ \geq \ {\tau \over \Delta(G)} \sum \Sb w:\ \{w, v\} \in E \\ l:\ 1\leq l \leq k \endSb
\left|z^{vw}_{il}\right| \left| {\partial \over \partial z^{vw}_{il}} \QQ^W_I(Z) \right|.$$
\medskip
(6.3.$r$) Let $W$ be an admissible sequence of vertices and let $I$ be an admissible sequence of indices such that 
$|W|=|I|=r$. Assuming that $W=(W', u, v)$ and $I=(I', j, i)$, we have
$$\left| {\QQ^{(W', u, v)}_{(I',j, i)}(Z) \over \QQ^{(W', u, v)}_{(I', i, j)}(Z)} -1 \right| \ \leq \ \epsilon.$$
\medskip
(6.4.$r$) Let $W$ be an admissible sequence of vertices and let $I$ be an admissible sequence of indices such that 
$|W|=|I|=r-1$. Let $w$ be a vertex such that the sequence $(W, w)$ is admissible and let $i$ and $j$ be indices such that 
the sequences $(I, i)$ and $(I, j)$ are admissible. Then $\QQ^{(W, w)}_{(I, i)}(Z) \ne 0$, $\QQ^{(W,w)}_{(I,j)}(Z) \ne 0$ and the angle between the complex numbers $\QQ^{(W, w)}_{(I, i)}(Z)$ and $\QQ^{(W,w)}_{(I,j)}(Z)$ does not exceed $\theta$.
\medskip
(6.5.$r$) Let $W$ be an admissible sequence of vertices and let $I$ be an admissible sequence of indices such that 
$|W|=|I|=r-1$. Let $u$ and $v$ be vertices such that the sequences $(W, u)$ and $(W, v)$ are admissible and let $i$  be an index such that 
the sequences $(I, i)$ is admissible. Then $\QQ^{(W, u)}_{(I, i)}(Z) \ne 0$, $\QQ^{(W,v)}_{(I,i)}(Z) \ne 0$ and the angle between the complex numbers $\QQ^{(W, u)}_{(I, i)}(Z)$ and $\QQ^{(W,v)}_{(I,i)}(Z)$ does not exceed $\theta$.
\bigskip
Let $W$ and $I$ be admissible sequences such that $|W|=|I|=r$, so that 
$W=\left(v_1, \ldots, v_r\right)$ and $I=\left(i_1, \ldots, i_r\right)$. If $r=|V|$ then 
$$\QQ^W_I(Z)=\prod \Sb 1 \leq j<l \leq r \\ \{v_j, v_l\} \in E \endSb z^{v_j v_l}_{i_j i_l} \ne 0,$$
so (6.1.$r$) holds for $r=|V|$. Moreover, denoting by $\deg\left(v_r\right)$ the degree of $v_r$, we obtain in this case  
$$ \sum \Sb w:\ \{w, v\} \in E \\ l:\ 1\leq l \leq k \endSb
\left|z^{vw}_{il}\right| \left| {\partial \over \partial z^{vw}_{il}} \QQ^W_I(Z) \right| =\deg\left(v_r\right) \left| \QQ^W_I(Z)\right|,$$
so (6.2.$r$) holds as well. Lemma 5.2 then implies that (6.3.$r$) holds for $r=|V|$, while (6.4.$r$) and (6.5.$r$) hold trivially, since for 
$W'=\left(v_1, \ldots, v_{r-1}\right)$ the only admissible extension is $W=\left(v_1, \ldots, v_r\right)$ and for $I'=\left(i_1, \ldots, i_{r-1}\right)$ the only admissible extension is $I=\left(i_1, \ldots, i_r\right)$. Hence (6.1.$r$)--(6.5.$r$) are satisfied for $r=|V|$.

From formula (3.2.1) and Lemma 4.1 we get the implication 
$$(6.1.r) \quad \text{and} \quad (6.4.r) \qquad \Longrightarrow \qquad  (6.1.r-1).$$
From Lemma 5.4 we get the implication 
$$ (6.4.r) \qquad \Longrightarrow \qquad  (6.2.r-1).$$
From Lemma 5.2 we get the implication
 $$ (6.1.r-1)  \quad \text{and} \quad (6.2.r-1) \qquad \Longrightarrow \qquad (6.3.r-1).$$
 From Part (1) of Proposition 4.3 we get the implication 
 $$(6.3.r) \quad \text{and} \quad  (6.4.r)  \qquad \Longrightarrow \qquad (6.5.r-1).$$
 From Part (2) of Proposition 4.3 we get the implication
 $$(6.3.r) \quad \text{and} \quad (6.5.r) \qquad \Longrightarrow \qquad (6.4.r-1).$$
This proves (6.1.2)--(6.5.2). From Part (2) of Proposition 4.3, we get the implication 
$$(6.3.2) \quad \text{and} \quad (6.5.2)  \Longrightarrow (6.4.1).$$
Thus for every vertex $v$ and every two indices $i$ and $j$ we have $\QQ^v_i(Z) \ne 0$, $\QQ^v_j(Z) \ne 0$ and the angle between complex numbers $\QQ^v_i(Z)$ and $\QQ^v_j(Z)$ does not exceed $\theta$. By formula (3.2.1) applied to the empty sequences $W$ and $I$ and Lemma 4.1, we conclude that $\QQ(Z) \ne 0$. 
{\hfill \hfill \hfill} \qed

\head Acknowledgments \endhead

The authors are grateful to Boris Bukh for suggesting Lemma 4.1 and to anonymous referees for careful reading of the paper, suggestions and corrections.

\enddocument

\end